\def\mytitle{A modified Bellman-Ford Algorithm for Application\texorpdfstring{\\}~in Symbolic Optimal Control and Plan and Goal Recognition}
\def\myname{Marcus Kreuzer, Alexander Weber and Alexander Knoll}
\def\mykeywords{Symbolic Optimal Control, Bellman-Ford Algorithm, Discrete Abstractions, High-Performance Computing, Plan and Goal Recognition}
\def\arxiv{0}
\documentclass[letterpaper, 10pt, conference]{ieeeconf}
\IEEEoverridecommandlockouts
\usepackage{cite}
\usepackage[english]{babel}
\usepackage{amsmath}
\usepackage{amssymb}

\usepackage{amsthm}
\usepackage{stmaryrd}
\usepackage{tikz}
\usetikzlibrary{calc,shapes,arrows,automata}
\usepackage{multirow}
\usepackage[]{caption} 
\usepackage{subcaption}
\usepackage{algorithm}
\usepackage{algpseudocode}
\algnewcommand\algorithmicinput{\textbf{Input:}}
\algnewcommand\Input{\item[\algorithmicinput]}
\swapnumbers
\newtheoremstyle{rem}{\topsep}{\topsep}{\normalfont}{0pt}{\bfseries}{.}{ }{\thmname{#1 }\thmnumber{#2}\thmnote{ \textup{(#3)}}}
\newtheorem{definition}{Definition}[section]

\theoremstyle{rem}
\newtheorem{example}[definition]{\pushQED{\qed}Example}
\makeatletter
\def\endexample{\popQED\@endtheorem}
\def\notion#1%
{\@nomath\begriff\ifdim\fontdimen\@ne\font>\z@%
\textbf{#1}\else\textit{#1}\fi}
\def\intcc#1{\ensuremath{\left[#1\right]}}

\makeatother
\def\powerset#1%
{\mathcal{P}(#1)%
}

\newcommand{\Frontier}{\mathcal{F}}
\def\implies{\relax\ifmmode\mathrel{\Rightarrow}\else$\implies$ \fi}
\usepackage{paralist}
\usepackage{hyperref}
\ifx\arxiv\undefined
\hypersetup{
colorlinks=true,%
breaklinks=true,%
pdfdisplaydoctitle=true,%
linkcolor={black},%
citecolor={black},
urlcolor={blue},
pdfstartview={FitV},%
pdftitle={\mytitle},
pdfauthor={\myname},%
pdfsubject={--},
pdfkeywords={\mykeywords}
}
\else
\hypersetup{
colorlinks=true,%
breaklinks=true,%
pdfdisplaydoctitle=true,%
linkcolor={black},%
citecolor={black},
urlcolor={blue},
pdfstartview={FitV},%
pdftitle={\mytitle},
pdfauthor={\myname},%
pdfsubject={Submitted to arxiv on \today.},
pdfkeywords={\mykeywords}
}
\fi
\ifx\arxiv\undefined
\else
\usepackage{fancyhdr}
\fi
\graphicspath{{figures/}} 
\title{\bf \LARGE \mytitle}
\author{
\myname
\thanks{
The authors are with the
Munich University of Applied Sciences,
Dept. of Mechanical, Automotive and Aeronautical Eng.,
80335 M\"unchen, Germany. 
}
\thanks{This work has been supported by the Bavarian Ministry of Economic Affairs, Regional Development and Energy (HAM-2307-0004).}%
}%
\begin{document}
\maketitle

\begin{abstract}
The contributions of this 
short technical note are two-fold.
Firstly, we introduce a modified version of 
a generalized Bellman-Ford algorithm 
calculating the value function
of optimal control problems 
defined on hyper-graphs. 
Those Bellman-Ford algorithms can be used 
in particular for the 
synthesis of near-optimal controllers 
by the principle of symbolic control. 
Our modification causes less nodes of
the hyper-graph being iterated during the 
execution compared to our initial 
version of the algorithm published in 2020. 
Our second contribution lies in the 
field of Plan recognition applied to
drone missions driven by symbolic controllers. 
We address and resolve 
the Plan and Goal Recognition monitor's dependence on 
a pre-defined initial guess 
for a drone's task allocation and mission execution.
To validate the enhanced implementation, 
we use a more challenging scenario for UAV-based aerial firefighting, 
demonstrating the practical applicability and robustness of the system architecture.
\end{abstract}
%
%
%
\section{Introduction}
\ifx\arxiv\undefined
\else
\thispagestyle{fancy}
\fi
Symbolic control synthesis refers to an 
automated method for computing controllers
for non-linear plants and quite general 
control objectives.
The synthesis process requires
firstly a so-called symbolic model 
of the plant, which may be interpreted 
as a discrete state machine. 
The synthesis on the controller 
is done based on this symbolic model and on 
a similarly transferred control objective. 
When this auxiliary synthesis is successful, 
the obtained controller can be refined 
to a controller for the actual control problem.
Roughly since twenty years, 
there is considerable progress in 
improving this technique, see e.g. 
\cite{ReissigWeberRungger17,WeberMacoveiciucReissig22}

In \cite{WeberKreuzerKnoll20} we introduced 
a generalized Bellman-Ford algorithm 
for the use in symbolic optimal control
motivated by the fact that in this way 
the synthesis step can be parallelized 
in contrast of algorithms based 
on Dijkstra's idea of a priority queue. 
In this note we present 
a modification of that algorithm
that reduces the number of processed nodes
in each iteration. 
As a second contribution, we show that our plan recognition (PR)
methodology introduced in 
\cite{KreuzerWeberLeupolzKnoll23,KreuzerWeberKnoll24}
can be successfully applied to other more challenging scenarios.
By strategically positioning target candidates in close proximity, 
we intentionally challenge the Plan- and Goal Recognition algorithm 
to operate at the limit of achievable accuracy \cite{Bitcraze_lighthouse_2024} 
with respect to the UAV system's position in space.

This technical note shall serve as a preliminary report
to a more extensive study in this area of research 
to be presented in near future.
The remaining part of the note is organized as follows.
The optimal control problem to consider is defined in Sect.~\ref{s:optimalcontrol}.
The modification of the algorithm in \cite{WeberKreuzerKnoll20} is presented in Sect.~\ref{s:algorithm}.
Then, experimental evaluations of the modified algorithm are presented (Sect.~\ref{s:experiments}). 
The last section of the paper (Sect.~\ref{s:PGRM}) deals with Plan- and Goal Recognition.
\section{The definition of Plant and Optimal Control Problem}
\label{s:optimalcontrol}
In this note, a \notion{plant} is a triple
\begin{equation}
\label{e:plant}
(X,U,F)
\end{equation}
of two non-empty sets $X$ and $U$, 
and a strict\footnote{A map $f$ with domain $A$ is \notion{strict} if $f(a) \neq \emptyset$ for all $a \in A$.} 
set-valued map 
$F\colon X \times U \rightrightarrows X$.
The time-discrete dynamics is given by 
\begin{equation}
\label{e:dynamics:discrete}
x(t+1) \in F(x(t),u(t))
\end{equation}
so that $X$ and $U$ are called 
the \notion{state} and \notion{input space}, 
respectively, of the plant.
In fact, for discrete sets $X$ and $U$ the triple \eqref{e:plant} defines a hyper-graph, 
whose transitions are determined via \eqref{e:dynamics:discrete}.

This note deals with
optimal control problems, 
which are defined by a 5-tuple
\begin{equation}
\label{e:def:ocp}
(X,U,F,G,g).
\end{equation} 
In \eqref{e:def:ocp}, the first entries are as in \eqref{e:plant} and the last
entries are maps
\begin{subequations}
\label{e:cost}
\begin{align}
G &\colon X \to \mathbb{R} \cup \{ \infty\}, \\
g &\colon X \times X \times U \to \mathbb{R} \cup \{\infty\},
\end{align}
\end{subequations}
which define the cost functional. 
Specifically, the objective is the one of an optimal stopping problem \cite{ReissigRungger13}. 
For time-discrete signals 
$x \colon \mathbb{Z}_+ \to X$, 
$u \colon \mathbb{Z}_+ \to U$ and 
$v \colon \mathbb{Z}_+ \to \{0,1\}$ 
the cost functional is defined by
\begin{equation}
\label{e:costs}
J(u,v,x) = G(x(T)) + \sum_{t = 0}^{T-1} g(x(t),x(t+1),u(t))
\end{equation}
if 
$T := \inf v^{-1}(1) < \infty$ and 
$J(u,v,x) = \infty$ if $v = 0$. 
Here, first edge from $0$ to $1$ in $v$
determines the stopping time $T$ 
of the plant.
Subsequently, the sets $X$ and $U$ in \eqref{e:def:ocp}
are assumed to be finite.
The algorithm to present computes the maximal fixed point
for the famous \notion{dynamic programming operator}
$P \colon \intcc{-\infty,\infty}^X \to \intcc{-\infty,\infty}^X$ defined by
\begin{equation}
\label{e:dynamicprogramming}
P(W)(x) = \min \left \{ G(x), Q(W,x) \right \}
\end{equation}
where 
\begin{equation}
\label{e:Q}
Q(W,x) := \inf_{u \in U} \sup_{y \in F(x,u)} g(x,y,u) + W(y).
\end{equation}
Here, the set $\intcc{-\infty,\infty}^X$ is the set of all maps $X \to \intcc{-\infty,\infty}$. 

In addition, for finite $X$ and $U$ a controller 
\begin{equation}
\label{e:controller}
\mu \colon X \to U
\end{equation} 
is \notion{optimal} if 
it realizes the minimum in \eqref{e:Q} when $W$ is the maximal fixed point.

Please refer to \cite{ReissigWeberRungger17,ReissigRungger18} for all the details 
to the previous definitions.
\section{The modified Algorithm}
\label{s:algorithm}
Algorithm \ref{alg:BellmanFord} is a modification of \cite[Alg.~1]{WeberKreuzerKnoll20}. 
(In the statement of Algorithm \ref{alg:BellmanFord} the set $\operatorname{pred}(x,u)$ is defined as
\begin{equation}
\label{e:pred}
\operatorname{pred}(x,u) = \{ y \in X \mid x \in F(y,u)\}
\end{equation}
which may be seen as the predecessors in the hyper-graph defined by $F$.)
In fact, lines 9--11, 16 and 19--20 are the novel extension. 
%
%
%
%
For reducing the number of processed nodes in each iteration, which is number of elements in $\mathcal{F}_1$, 
the modification used the fact that a node $x$ may obtain a new value $W(x)$ if the supremum 
of its successors changes. Therefore, the supremum is stored in line 16 and the upcoming frontier is extended only if the previously explain condition is true (line 19). 
\begin{algorithm}
\caption{\label{alg:BellmanFord}\small{Modification of the Bellman-Ford Algorithm in \cite{WeberKreuzerKnoll20}}}
\begin{algorithmic}[1]
\Input{Optimal control problem $(X,U,F,G,g)$}
\State{$\Frontier_1 \gets \emptyset$\hfill{}{// ``Active" Frontier \cite{CormenLeisersonRivestStein09}}}
\State{$\Frontier_2 \gets \emptyset$\hfill{}{// ``Upcoming" Frontier \cite{CormenLeisersonRivestStein09}}}
\ForAll{$x \in X$}
\State{\label{alg:BF:init:V}$W(x) \gets G(x)$}
\If{$G(x) < \infty$}
\State{$\Frontier_1 \gets \Frontier_1 \cup ( \cup_{u \in U } \operatorname{pred}(x,u) )$}
\EndIf{}
\EndFor{\label{alg:init:endfor}}
\ForAll{$(x, u) \in X \times U$}
\State{$F_{\max}(x,u) \gets \{y\}$ for some $y \in F(x,u)$}
\EndFor{}
\State{$i=0$}
\While{\label{alg:BF:while}$\Frontier_1 \neq \emptyset$ \textbf{and} $i < |X|$}
\For {\label{alg:BF:for}$(x,u) \in \Frontier_1 \times U$}
\State{\label{alg:BF:sup}$d \gets \sup_{y \in F(x,u)} ( \ g(x,y,u) + W(y) \ ) $}
\State{\label{alg:BF:storemax}\textcolor{black}{$F_{\max}(x,u) \gets \operatorname{arg\,sup}\{ W(y) \mid y \in F(x,u)\}$}}
\If{$d < W(x)$}
\State{\label{alg:BF:assign}$W(x) \gets d$}
\State{\label{alg:BF:conditionnew:1}$\mathcal{F}_{\max} \gets \{ y \in X \mid \exists_{\tilde u \in U} \ x \in F_{\max}(y,\tilde u) \}$}
\State{\label{alg:BF:conditionnew:2}\textcolor{black}{$\Frontier_2 \gets \Frontier_2 \cup \mathcal{F}_{\max}$}}
\EndIf{}
\EndFor{\label{alg:BF:endfor}}
\State {$\Frontier_1 \gets \Frontier_2$\hfill{}{// Swap frontiers}}
\State{$\Frontier_2 \gets \emptyset $}
\State{$i \gets i + 1$}
\EndWhile{}
\State{\Return{$W$}}
\end{algorithmic}
\end{algorithm}
\section{Experimental evaluation}
\label{s:experiments}
Algorithm \ref{alg:BellmanFord} was tested of three examples for this note, whose 
results are summarized in Tab.~\ref{tab:keyfigures}. 
The focus of the analysis is intentionally 
not on run-times and implementations but on qualitative understanding.
As seen from Tab.~\ref{tab:keyfigures} the number of processed nodes is in every problem significantly reduced.
However,
the total number of required iteration increases. 
With our current implementation in software, run-times remain almost identical to the ones using the initial version of the algorithm. 
Therefore, future research will focus on increasing the efficiency of a single iteration. ´
\begin{table}
\centering
\begin{tabular}{cccccc}
Reference & Problem no. & $i$ & $\sum_i \mathcal{F}_1 / | X | $ \\
\hline 
\multirow{2}{7em}{\cite[Sect.~IV.A]{WeberKnoll20}} 
& 1 & 99 (86) & 8.97 (17.01)\\
& 2 & 85 (81) & 8.98 (13.23) \\
 \hline
\multirow{5}{7em}{\cite{KreuzerWeberLeupolzKnoll23}}
& 1 & 108 (77) & 8.71 (11.38)\\
& 2 & 92 (53) & 7.37 (9.23)\\
& 3 & 118 (111) & 8.72 (18.38)\\
& 4 & 94 (60) & 8.35 (11.14)\\
& 5 & 104 (57) & 8.62 (11.27)\\
\hline
\multirow{5}{7em}{Example in Sect.~\ref{s:PGRM} of this paper} 
& 1 & 89 (58) & 6.28 (6.45)\\
& 2 & 70 (45) & 4.25 (5.65)\\
& 3 & 78 (44) & 4.79 (5.44)\\
& 4 & 78 (43) & 5.06 (5.80)\\
& 5 & 44 (45) & 3.09 (3.47)\\
\end{tabular}
\caption{\label{tab:keyfigures}Qualitative key figures of Alg. \ref{alg:BellmanFord} in comparison to \cite[Alg. I]{WeberKreuzerKnoll20}.}
\end{table}
\section{Application to Mission Monitoring}
\label{s:PGRM}
Building upon the mission control framework and 
the Plan and Goal Recognition Monitor (PGRM) developed for 
single \cite{KreuzerWeberLeupolzKnoll23} and 
multi-UAV systems \cite{KreuzerWeberKnoll24}, 
this work introduces the following key advancement: 
The need for an initial guess for task assignment 
concerning mission execution is eliminated. 
Instead, the PGRM dynamically calculates 
the first estimated sequence of mission subtasks 
for an observed UAV system following its takeoff phase.
Furthermore, through this case study, 
we demonstrate the general applicability of the presented methods 
in the aerial domain by validating their effectiveness 
in a scenario characterized by closely spaced targets 
(compared to 
\cite{KreuzerWeberLeupolzKnoll23,
KreuzerWeberKnoll24}), 
again considering spatio-temporal 
uncertainty inherent to cyber-physical systems.\\
This section provides a brief overview 
of the relevant topics and methodologies. 
Subsequently, the scenario is outlined, 
and the mission's progression 
is detailed along its timeline.\\
\subsection{Plan Recognition for prediction of trajectory and mission (sub-)tasks} 
The concept of \textit{plan recognition by planning}, 
introduced in 
\cite{ramirez2009plan,
ramirez2010probabilistic}, 
enables an observer to infer a system's behavior 
by evaluating its actions relative to a set of virtual, 
hypothetical plans produced 
by a rational (cost-based), deterministic planning algorithm.
Inspired by aforementioned works, 
\cite{FitzpatrickLipovetzkyPapasimeonRamirezVered21,
KaminkaVeredAgmon18,
masters2019cost,
MastersSardina18,
vered2016online} extended the method to 
the navigational domain based on 
RRT \cite{sucan2012open} and A* \cite{pohl1970heuristic} algorithms.\\ 
By combination with the symbolic optimal control (SOC) framework, 
we made the method applicable to continuous surveillance 
for airspace coordination \cite{KreuzerWeberLeupolzKnoll23} and 
decentralized Runtime Assurance for mission execution \cite{KreuzerWeberKnoll24}. 
In this context, for plan recognition 
the PGRM evaluates the alignment between calculated trajectory candidates 
(each leading to an available mission goal) and a UAV's online, 
incrementally observed partial trajectory.
The SOC framework serves two primary functions:\\
Firstly, it operates the closed loop of a controller $\mu_k$ to guide an UAV system along its trajectory toward one of $|K|$ mission (sub-)targets, also referred to as \textit{goals}.\\
Secondly, the SOC algorithm facilitates a model-based, \textit{faster-than-real-time} online simulation. This procedure virtually executes the controllers $\mu_k$ on an UAV system's model to generate the aforementioned set of $|K|$ hypothetical trajectory candidates, which are subsequently analyzed within the plan and goal recognition process.\\     
We briefly outline the problem based on our work in 
\cite{KreuzerWeberLeupolzKnoll23,
KreuzerWeberKnoll24}:\\  
Assume a strict set-valued map 
$H \colon X \rightrightarrows Y$ 
exists, mapping the states of the plant \eqref{e:plant} to subsets of a specified non-empty output set $Y$. 
Additionally, consider a potentially incomplete observation sequence 
$(o(0),\ldots,o(T))$ 
within $Y$, satisfying the following property:
For every $t \geq 0$, $t \leq T$, 
there exists an initial state $p \in X$, a controller $\mu$ as specified in \eqref{e:controller}, and a pair $(u,x)$
of signals $x \colon \mathbb{Z}_+ \to X$, $u \colon \mathbb{Z}_+ \to U$ such that
\begin{subequations}
\label{e:planrecognition:condition}
\begin{align} 
    x(t + 1) &\in F(x(t), \mu(x(t))), \\
	o(t) &\in H(x(t))
\end{align}
\end{subequations}
for all $t \leq T$ and some $T \in \mathbb{Z}_+$. 
Using this, the Plan Recognition (PR) problem can be formulated:
Given $K > 1$ controllers $\mu_1, \ldots, \mu_K$ for the plant \eqref{e:plant}, 
along with an observation sequence $(o(0),\ldots,o(T))$ for $T > 0$, 
identify the controller $\mu_{\ast} \in \{\mu_1, \ldots, \mu_K\}$ 
such that there exists $p \in X$ and a pair of signals $(u, x)$ as above 
satisfying \eqref{e:planrecognition:condition} 
for all $t \leq T$ with $\mu_\ast$ in place of $\mu$ in \eqref{e:planrecognition:condition}.

We address uncertainties and inaccuracies we experienced by real-world experiments (see \cite{KreuzerWeberLeupolzKnoll23}):
For the concrete example below, 
we assume $Y = \mathbb{R}^3$ and define $H$ as:
${H(x) = x + \intcc{-\kappa,\kappa}^3}$
where $\kappa = 0.075$ m, 
representing the diagonal half-span of the UAV. 
This setup accounts for, 
as well as simulates acceptable positional offsets.\\
Using $\kappa$ as a threshold parameter for trajectory course violation, we assume the real-world systems median Euclidean error $\epsilon < 0.04m$ \cite{TaffanelEtAl21}, regarding positional accuracy, as being considered. \\ 
By continuously monitoring the UAV's position over time and in case of a threshold violation,  
solving the PR problem during mission execution the PGRM identifies the current controller $\mu_{\ast}$. 
After acquiring this information, 
the PGRM is able to predict the future trajectory and, consequently, its ending state (=target), which is associated with an observed, 
cooperative UAV's current mission (sub-)goal.\\
Next, an UAV firefighting mission resulting in 
successful completion is presented 
in an experimental scenario demonstrating 
the previously described implementations. 
A timeline, 
qualitatively illustrating mission execution, 
is provided in Fig.\ref{fig:timeline_mission}.
\subsection{Multi-UAV firefighting scenario}
\begin{figure}
	\centering
	\includegraphics[scale=.53]{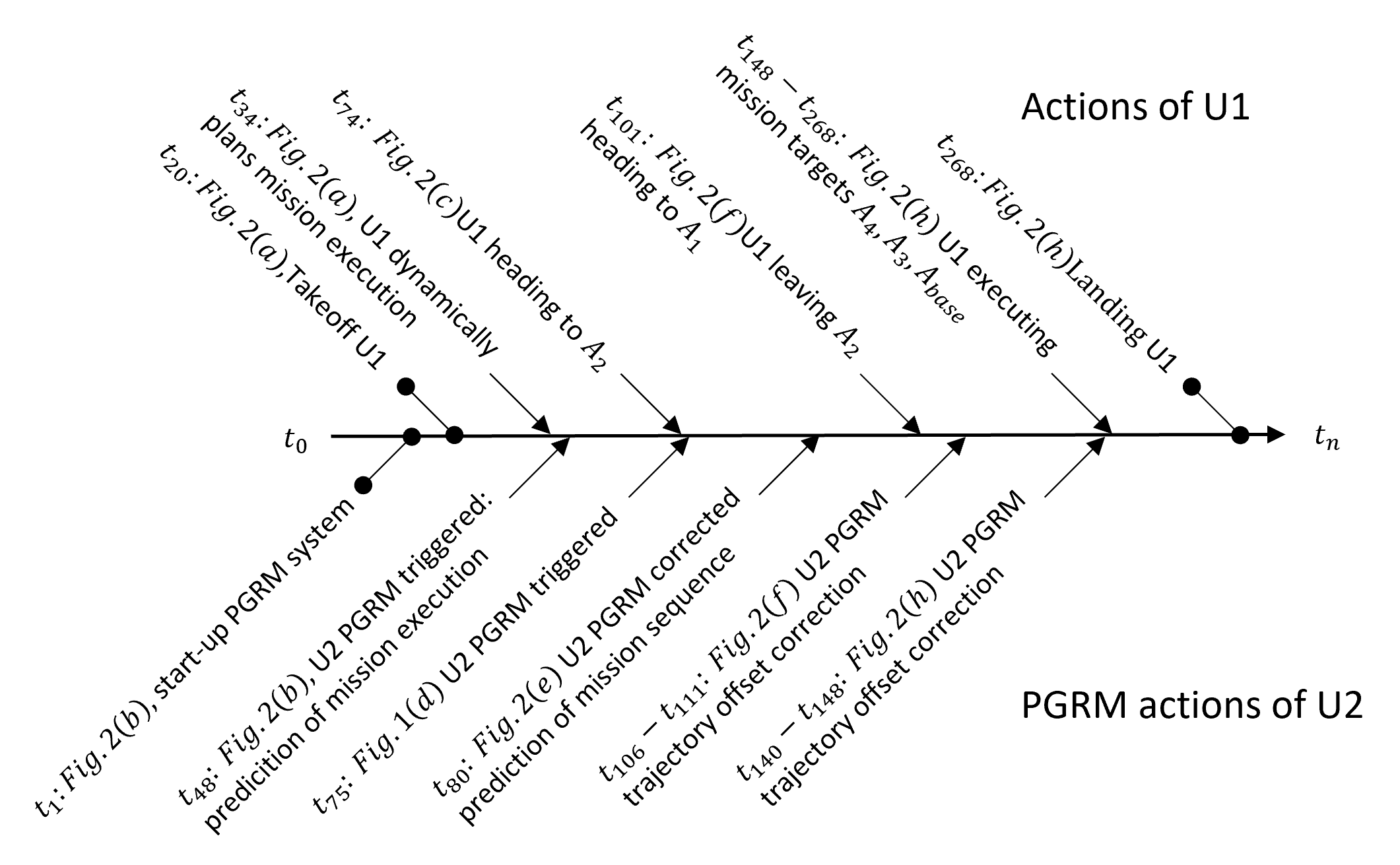}
	\caption{\label{fig:timeline_mission}Timeline for mission execution of U1 and the PGRM's processes to monitor correct trajectory alignment and (sub-)task assignment. Qualitative illustration of sequences including reference to figures.}
	\vspace*{-\baselineskip}
\end{figure}

\begin{figure*}
	\centering
	\begin{subfigure}[c]{0.49\textwidth}\hspace*{-1.1cm}
		\includegraphics[scale=.325 ]{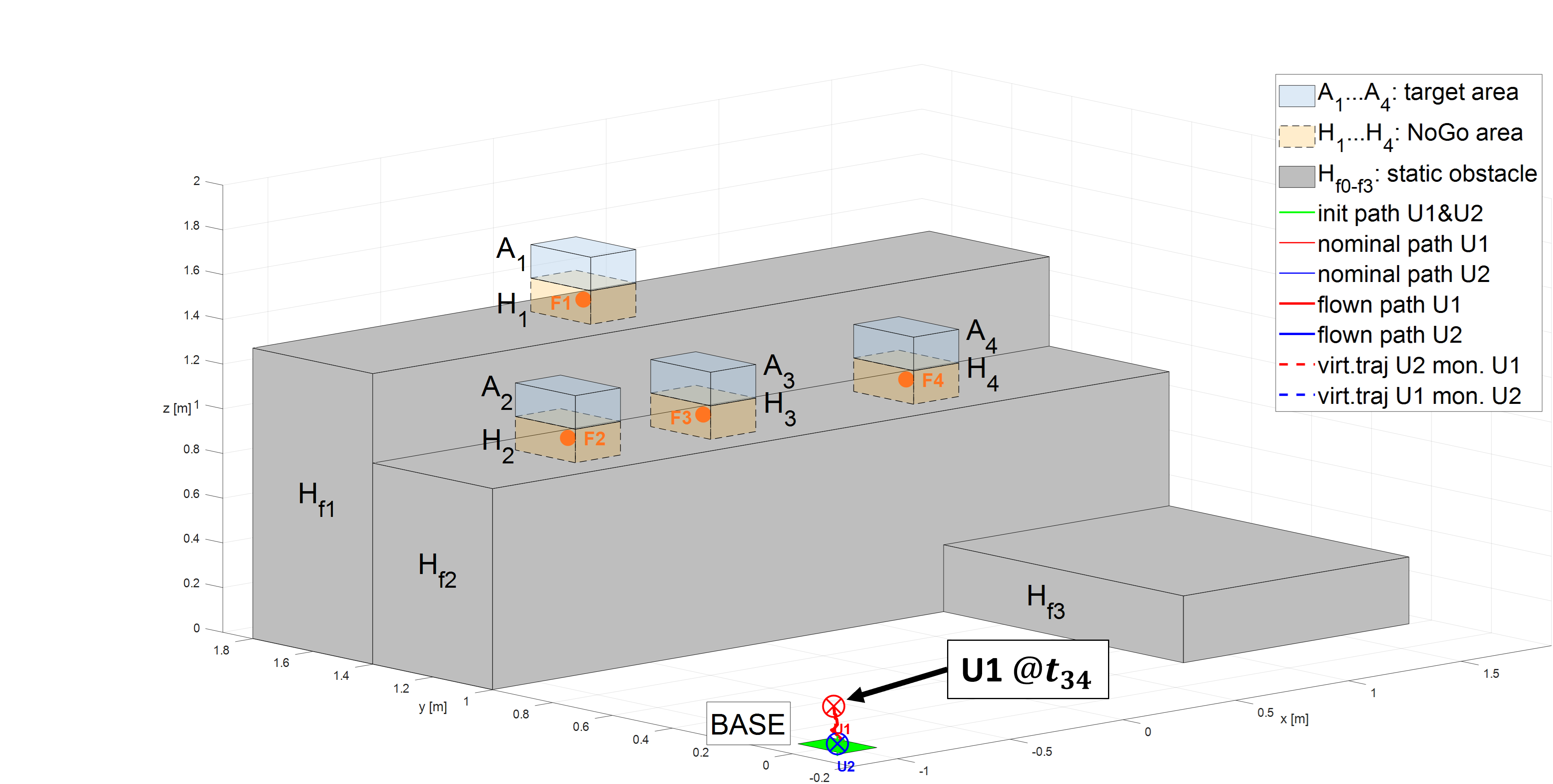} 
		\subcaption{\label{fig:U1_TO_done}Finished takeoff phase and dynamic mission planning of U1(red) $(@t_{34})$. For target sequence, $A_2\rightarrow A_1\rightarrow A_4\rightarrow A_3\rightarrow A_{base}$ is obtained. U2(blue), staying on ground, is monitoring U1 and thus ready for a plan recognition (PR) process.}
	\end{subfigure}\hfill
	\begin{subfigure}[c]{0.49\textwidth}\hspace*{-0.1cm}
		\includegraphics[scale=.3 ]{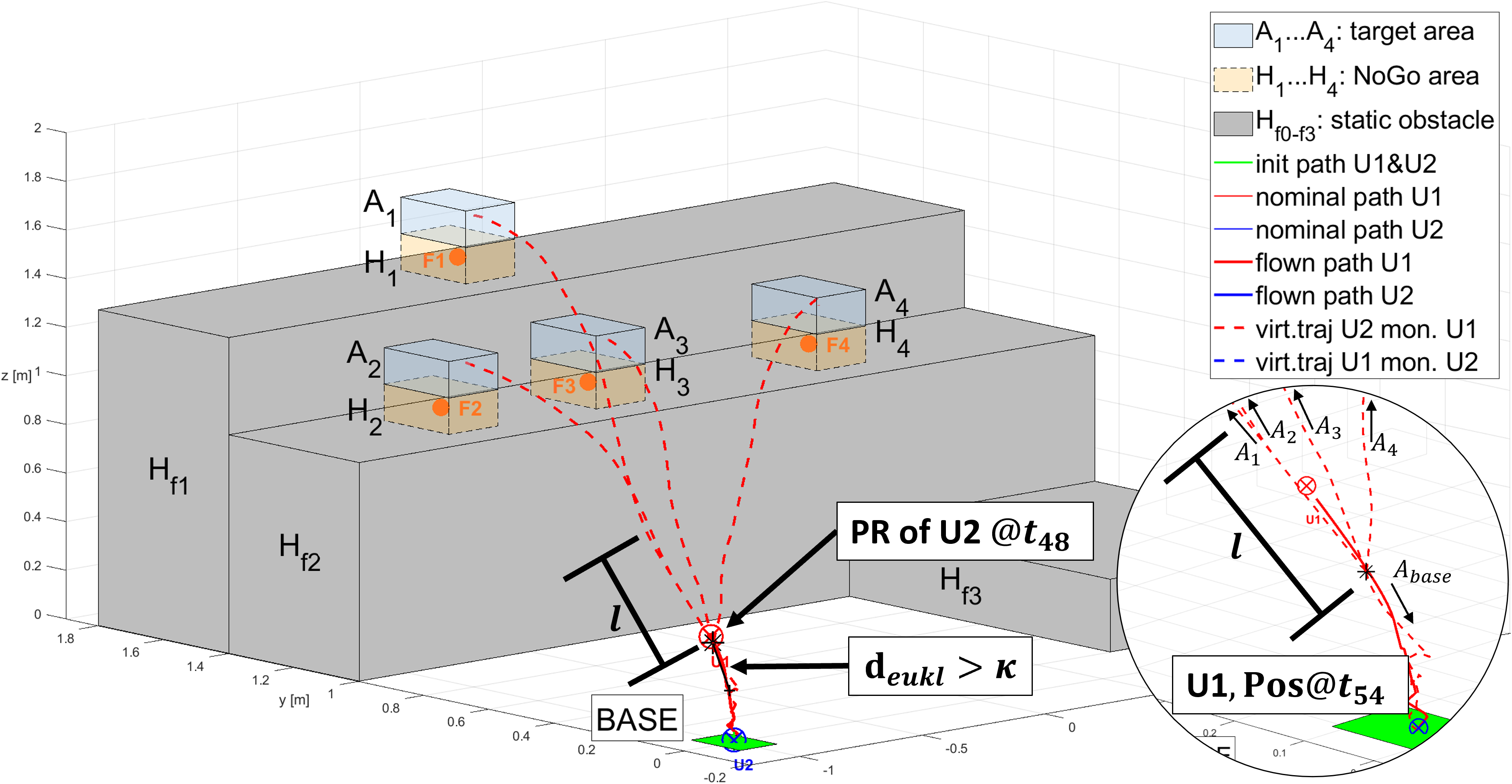}
		\subcaption{\label{fig:U2monU1_hypos_t2}U2's PGRM triggered by $d_{eukl.} > |\kappa|$, generating hypothetical trajectories for U1's target candidates (red dashed lines, black asterisk $@t48$), showing close alignment in the area of $l$ for trajectories to $A_1$ and $A_2$.  The zoomed area depicts U1 $@t_{54}$, following both near to identical tracks.}
	\end{subfigure}\hfill
	\centering
	\begin{subfigure}[c]{0.49\textwidth}\hspace*{-0.1cm}
		\includegraphics[scale=.3]{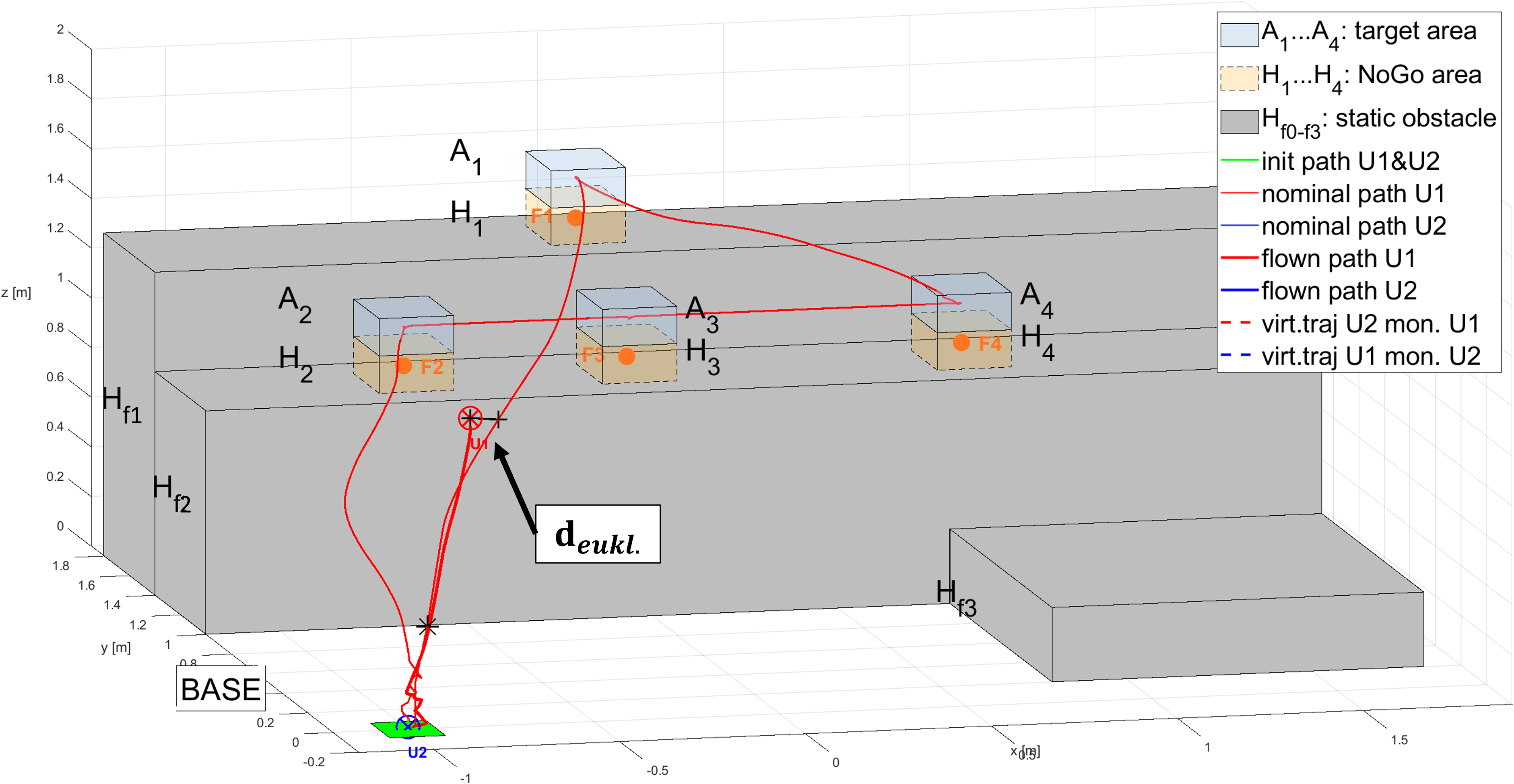}
		\subcaption{\label{fig:U1_hdg_to_A2_t4} $t_{74}$: U1 leaving near to identical trajectory courses heading to $A_2$. The PGRM'S incorrect identification of $A_1$ as being the target of U1 results in the growing distance $d_{eukl.}$.}
	\end{subfigure}\hfill
	\begin{subfigure}[c]{0.49\textwidth}\hspace*{-0.1cm}
		\includegraphics[scale=.3]{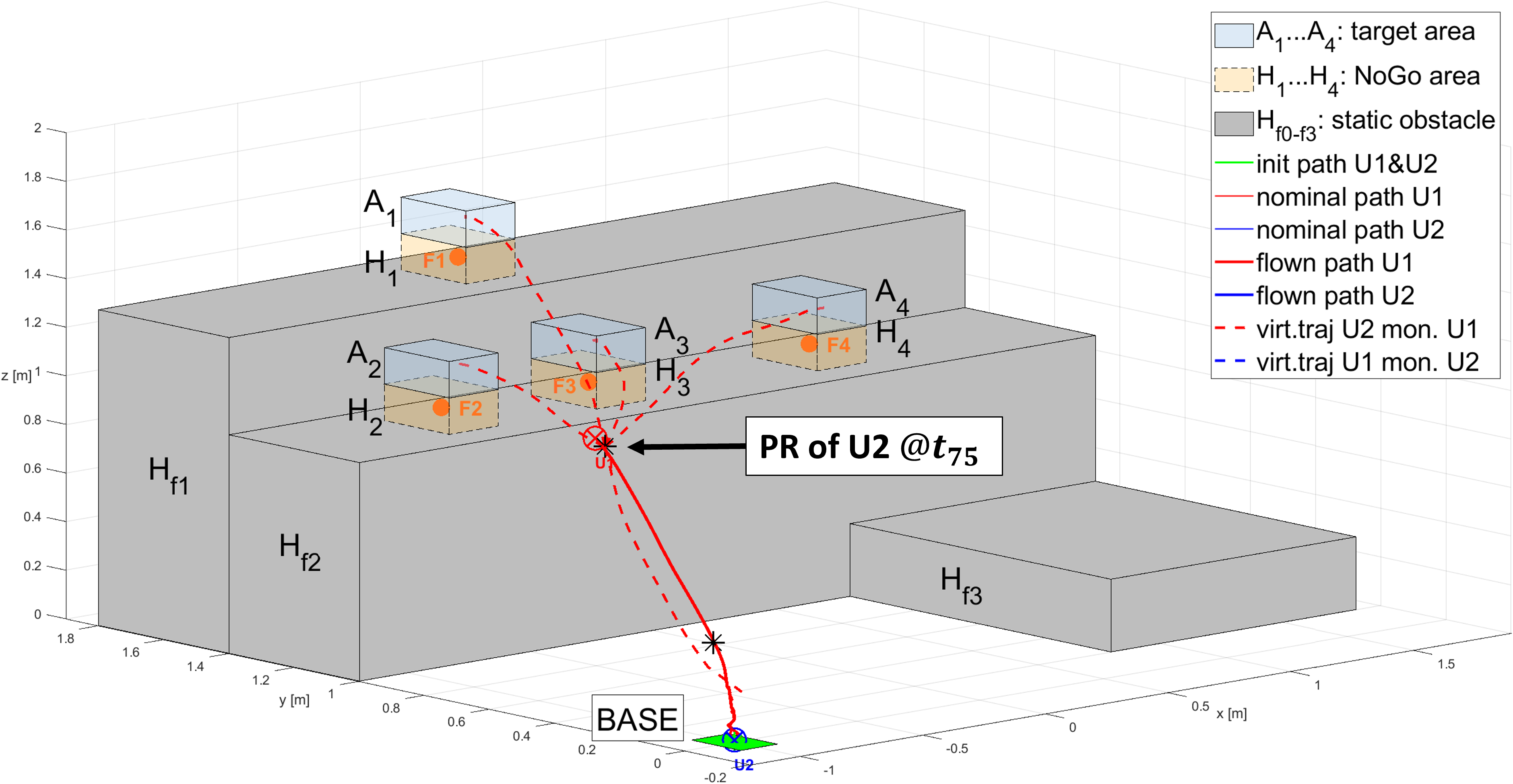}
		\subcaption{\label{fig:U1_hdg_to_A2_t5}In this situation $@t_{75}$, the hypothetical trajectory candidates show clear distinction in regards to their course to targets $A_1,...,A_4$ and $A_{base}$, see red dashed lines.}
	\end{subfigure}\hfill
	\centering
	\begin{subfigure}[c]{0.49\textwidth}\hspace*{-0.1cm}
		\includegraphics[scale=.3]{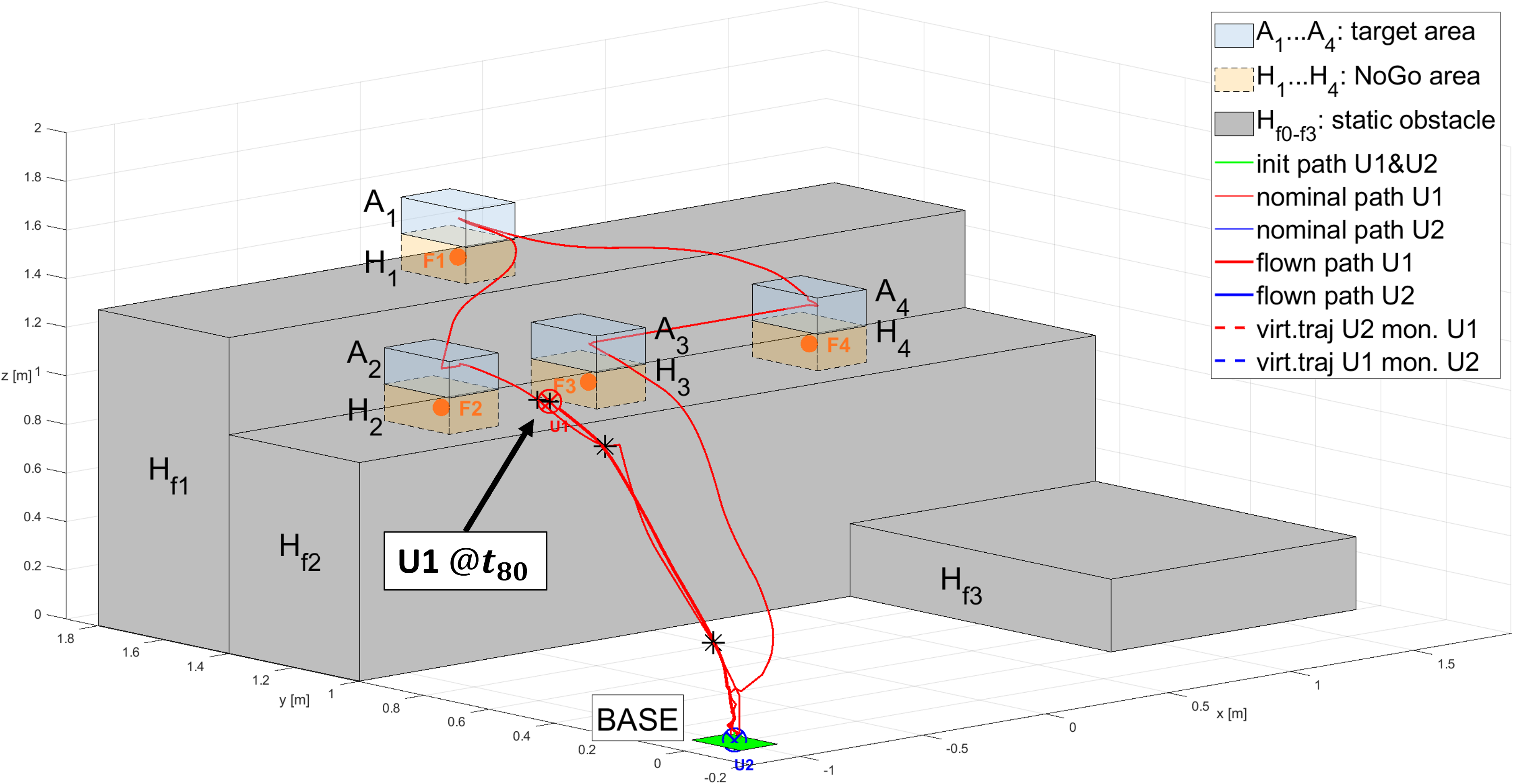}
		\subcaption{\label{fig:U1_hdg_to_A2_PR_done_t6} The PGRM of U2 infers the correct target $A_2$ of U1 ($@t_{80}$) with a slight acceptable offset in trajectory prediction, see black line. Notice $d_{eukl.}<\kappa$. Now, target sequence is predicted as $A_2\rightarrow A_1\rightarrow A_4\rightarrow A_3\rightarrow A_{base}$. }
	\end{subfigure}\hfill
	\begin{subfigure}[c]{0.49\textwidth}\hspace*{-0.1cm}
		\includegraphics[scale=.3]{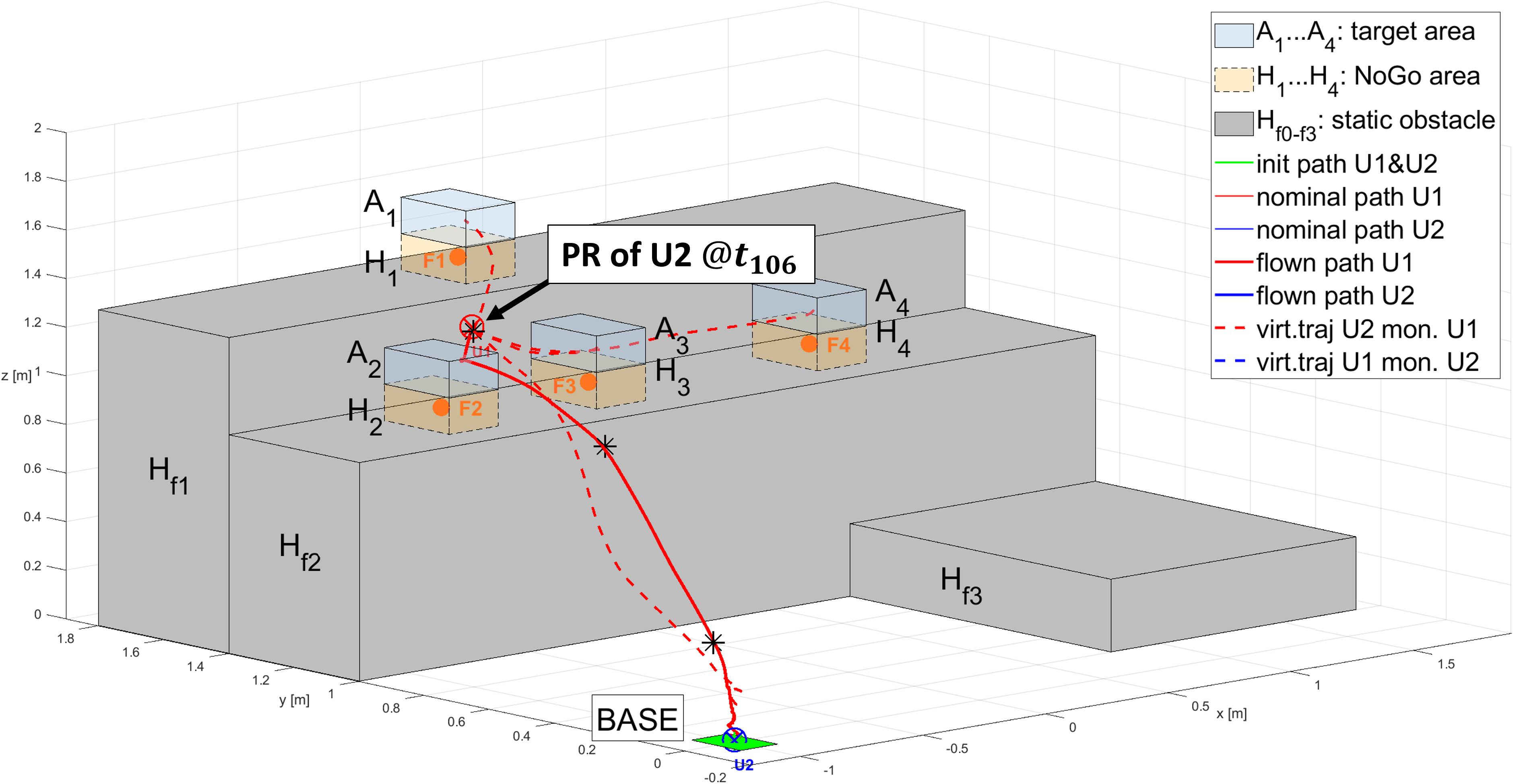}
		\subcaption{\label{fig:U1_hdg_to_A1_PR_offset_t7}After leaving $A_2$, U1 follows a track slightly offset from prediction, which triggers a PGRM process $@t_{106}$ to refine the trajectory estimation and assure $d_{eukl.}<\kappa$ with unchanged target $A_1$. Refer also to Fig.\ref{fig:U1_hdg_to_A1_PR_offset_done_t8}.} 
	\end{subfigure}\hfill
	\centering
	\begin{subfigure}[c]{0.49\textwidth}\hspace*{-0.1cm}
		\includegraphics[scale=.3]{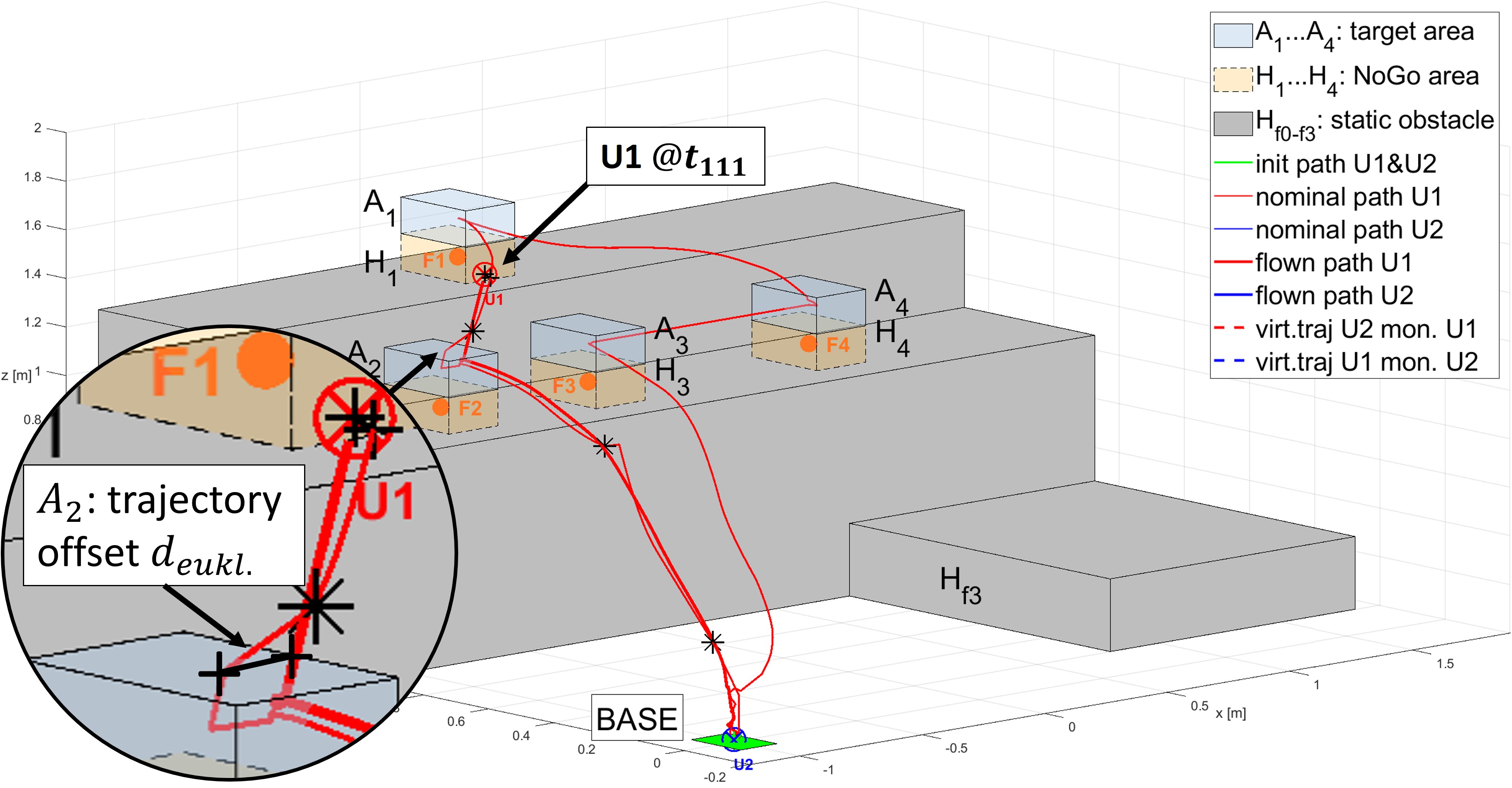}
		\subcaption{\label{fig:U1_hdg_to_A1_PR_offset_done_t8} The PGRM of U2 infers the correct target and trajectory of U1 $@t_{111}$. The zoomed area depicts the corrected trajectory offset, which initially triggered this PGRM process $@t_{106}$ (black asterisk).}
	\end{subfigure}\hfill
	\begin{subfigure}[c]{0.49\textwidth}\hspace*{-0.1cm}
		\includegraphics[scale=.3]{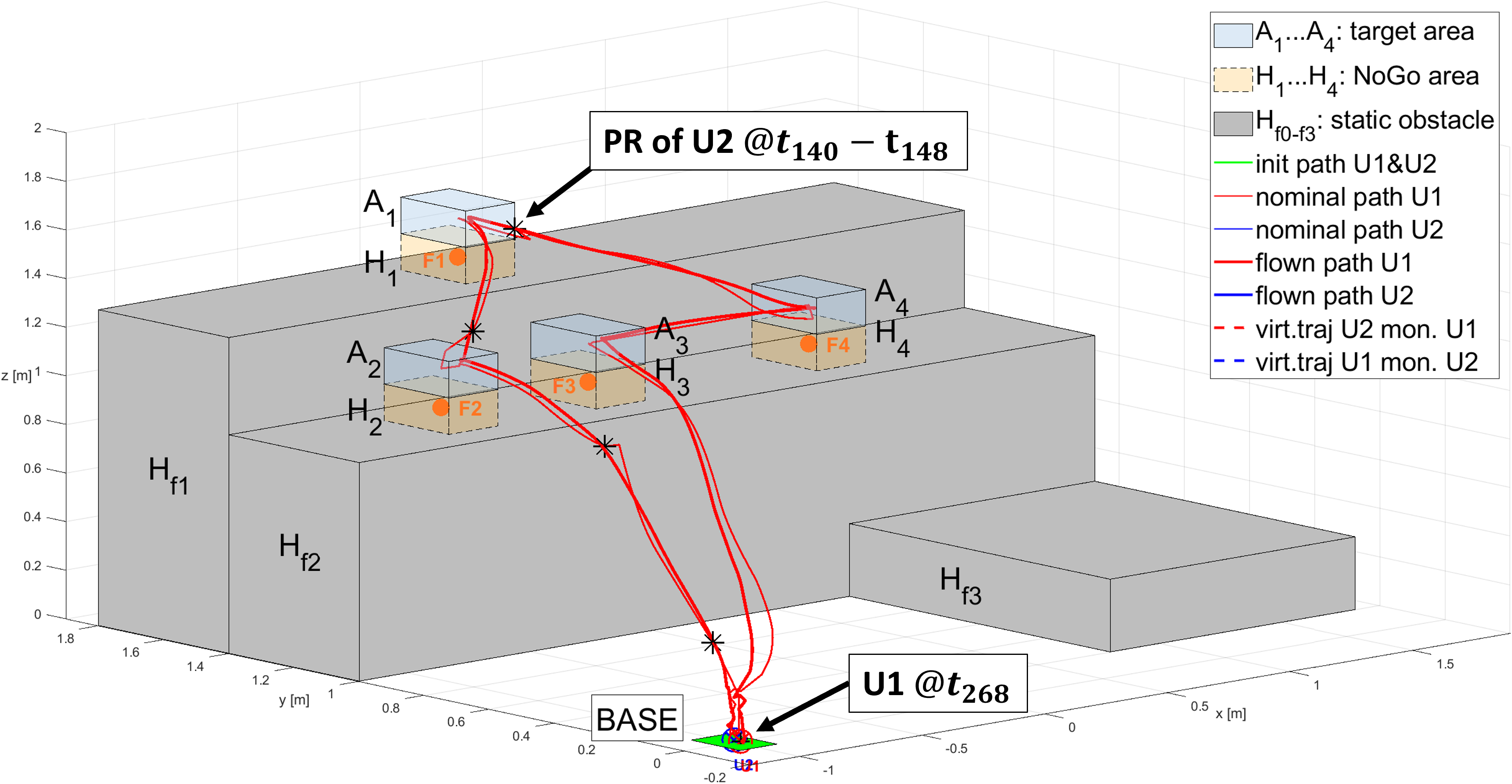}
		\subcaption{\label{fig:U1_touchdown_t268}U2's PGRM mitigates a second trajectory offset from $t_{140}$ to $t_{148}$. U1 finishes its mission by landing at the base $@t_{268}$. This scene represents the end state of successful mission execution.}
	\end{subfigure}\hfill
	\caption{\label{fig:crazy_mission}Mission execution shown in scenes.}
\end{figure*}

The mission's start- and end-state is at the base $A_{base}$. 
It is set within a scenario map featuring fixed obstacles $H_{f0},...,H_{f3}$ representing a mountainside, 
No-Go Areas (fires) $H_0,...,H_4$, 
and corresponding target areas $A_1,...,A_4$.
For successful mission completion, 
UAV1 (U1, red) must visit all target areas. 
See Tab.\ref{tab:basicproblemdata} and Fig.\ref{fig:crazy_mission}.
\begin{table}[H]
	\centering
	\begin{tabular}{|l|l|}
		\hline
		Symbol \& Meaning & Value (north-east-up coordinate system) \\
		\hline \hline
		$X_0$ (mission area) & $\intcc{-1.2,1.8}\!\times\!\intcc{-0.2,1.8}\!\times\!\intcc{-0.2,1.8}$ \\
		$A_{\textrm{base}}$ (base) & $\intcc{-1.1,-0.9}\!\times\!\intcc{-0.1,0.1}\!\times\!\{0\}$\\
		$A_{\textrm{hover}}$ (generic hover area) & $\intcc{-0.1,0.1}\!\times\!\intcc{-0.1,0.1}\!\times\!\intcc{0,0.225}$ \\
		$A_1$ (water release area) & $(0.0,1.65,1.525) + A_{\textrm{hover}}$\\
		$A_2$ (water release area) & $(0.6,1.2,1.125) + A_{\textrm{hover}}$\\
		$A_3$ (water release area) & $(0.0,1.2,1.125) + A_{\textrm{hover}}$\\
		$A_4$ (water release area) & $(0.9,1.2,1.125) + A_{\textrm{hover}}$\\
		$A_5$ (water release area) & $(1.3,0.6,0.525) + A_{\textrm{hover}}$\\
		$H_1,\ldots,H_5$ (fires)   & $ - (0,0,0.15) + A_i$, $i \in \{1,\ldots,5\}$\\
		$H_{f1}$ (fixed obstacle/hill)   & $\intcc{-1.2,1.8} \times \ldots$ \\
		
		& \phantom{\ldots}$\intcc{1.4,1.8} \times \intcc{-0.2,1.3}$\\
		$H_{f2}$ (fixed obstacle/hill)   & $\intcc{-1.2,1.8} \times \ldots$ \\
		
		& \phantom{\ldots}$\intcc{1.0,1.4} \times \intcc{-0.2,0.9}$\\
		$H_{f2}$ (fixed obstacle/hill)   & $\intcc{0.8,1.8} \times \ldots$ \\
		
		& \phantom{\ldots}$\intcc{0.2,1.0} \times \intcc{-0.2,0.3}$\\
		\hline
	\end{tabular}
	\caption{\label{tab:basicproblemdata}Spatial object data from the mission. Negative heights in $X_0$ result from the positioning system employed in real-world experiments, which produces a slightly tilted coordinate system, see \cite{KreuzerWeberLeupolzKnoll23}. Units are meters.}
	\vspace*{-\baselineskip}
\end{table}
The scenario in this work aims at 
testing the robustness of the plan recognition method 
against false positively identified trajectories, 
respectively mission tasks.\\
Unlike the scenarios presented in 
\cite{KreuzerWeberLeupolzKnoll23, 
KreuzerWeberKnoll24}, 
the target areas in this experiment are more closely spaced. 
In this way, 
for a PR process, 
we intentionally provoke the generation of a set of trajectory hypotheses 
characterized by a spatial proximity 
in the order of magnitude of $|\frac{\kappa}{2}|$, 
reflecting the system's threshold for positional uncertainty. \\
Within the recognition process, 
this will complicate the correct assignment of 
the currently observed UAV position to the matching hypothetical trajectory, 
which is mainly based on the evaluation of Euclidean distance.\\
We demonstrate that this issue is inherently mitigated through 
the PGRM's architecture within acceptable performance degradation: 
Influence on airspace management and 
coordination is exhibited for a short time period and 
thus is negligible in relation to safety concerns. \\ 
%
We now describe mission progress and 
actions taken by U1 while being monitored by the PGRM mechanism of UAV2 (U2, blue), 
which stays on ground. 
We follow the discretized timeline by time steps $t_n$:\\
After take-off at $t_{20}$ and 
the completion of an initial ascending phase at $t_{34}$, 
U1 dynamically plans the sequence of available targets to be reached 
(fires $F_1,...,F_4$, respectively target areas $A_1,...,A_4$), 
using the planning and execution architecture based on symbolic control as presented above and in \cite{KreuzerWeberLeupolzKnoll23, KreuzerWeberKnoll24}. 
See Fig.\ref{fig:U1_TO_done} and 
the timeline in Fig.\ref{fig:timeline_mission}. \\
As a result, 
U1 heads on to $A_2$, 
the first target of the correctly planned sequence 
$A_2\rightarrow A_1\rightarrow A_4\rightarrow A_3\rightarrow A_{base}$.\\
The PGRM subsystem of U2 monitors U1's trajectory course. 
Due to dynamic planning, 
there is no initial guess for 
the first target of U1.
The ``last" nominal state 
we know in advance for U1's trajectory 
is the final state of its ascending phase after take-off (see Fig.\ref{fig:U1_TO_done}, position of U1). 
At $t_{48}$ (see black asterisk in Fig.\ref{fig:U2monU1_hypos_t2}), 
U2 triggers a PGRM process, since the threshold $\kappa$ is exceeded 
by the Euclidean distance $d_{eukl.}$ (black line) between 
the aforementioned nominal and the observed trajectory state of U1. \\
By the PGRM, 
a set of hypothetical trajectories is generated to targets $A_1,...,A_4$ and $A_{base}$. 
We note that the proximity along the line segment $l$ of trajectories 
heading to $A_1$ and $A_2$ is below $\frac{\kappa}{2}$, 
as intended for this experiment. 
See Fig.\ref{fig:U2monU1_hypos_t2}. \\   
In the following time steps, 
U2's PGRM evaluates the alignment of U1's observed states and 
the states of hypothetical trajectories. 
Considering uncertainty of positional data inherent to the system, 
trajectories to $A_1$ and $A_2$ appear near to identical and 
the PGRM infers at $t_{54}$, 
that U1 is heading to $A_1$, 
which will turn out to be a faulty result. 
See zoomed area in Fig. \ref{fig:U2monU1_hypos_t2}.\\
Additionally to the current upcoming target, 
the PGRM predicts the future sequence of targets to be visited by U1 based on 
the solution of a travelling salesman problem, 
see \cite{KreuzerWeberLeupolzKnoll23} for details. 
This is illustrated by the thin red line touching areas 
$A_1\rightarrow A_4\rightarrow A_3\rightarrow A_2\rightarrow A_{base}$ in series. 
It will be monitored by the PGRM from now on.\\
From the perspective of U2's demands on safe airspace coordination, e.g. to keep separation minima,
this trajectory is representing the known, 
yet \textit{safe set} of U1’s mission trajectory (considering $\kappa$, too), 
see Fig.\ref{fig:U1_hdg_to_A2_t4}. \\
At $t_{74}$, the Euclidean distance between 
trajectory prediction to $A_1$ and U1's observed track 
(black line, state positions marked by "+") is about to 
exceed the threshold $\kappa$. 
Remember, U2's PGRM incorrectly inferred $A_1$ as U1's target instead of $A_2$. \\
Thus, at $t_{75}$, $d_{eukl.} > \kappa$ holds true and 
another recognition process is triggered by the PGRM of U2. 
In this situation, the generated trajectory candidates are characterized by a considerable margin alongside their course, see Fig.\ref{fig:U1_hdg_to_A2_t5}. \\
After $A_2$ was correctly identified as being the current target $@t_{80}$, 
the prediction for the sequential mission execution results in $A_2\rightarrow A_1\rightarrow A_4\rightarrow A_3\rightarrow A_{base}$. 
This is depicted by the thin red line, see Fig.\ref{fig:U1_hdg_to_A2_PR_done_t6}.\\
On U1's arrival at target area $A_2$, 
we simulate another real-world requirement: 
Assume, that U1 performs local, 
real-time surveillance of fire hot spots. \\
Typically, hot spots to be extinguished continue to move its position 
while eating their way through the landscape (also during mission time). 
Thus, the exact target position has to be corrected dynamically within the initially assumed target area, 
which then results in a distinct track for U1 to leave this area. 
Regarding real-world scenarios, this is expected to be the case regularly.\\
Thus in our experiment, 
heading to $A_1$, U1 leaves target area $A_2$ following a track slightly offset from prediction (see zoomed area in Fig. \ref{fig:U1_hdg_to_A1_PR_offset_done_t8}) and 
the PGRM of U2 is triggered $@t_{106}$, 
producing hypothetical trajectories 
to remaining targets $A_1, A_3, A_4, A_{base}$, see Fig.\ref{fig:U1_hdg_to_A1_PR_offset_t7}. \\
After completion of the PGRM process $@t_{111}$, 
the offset of the PGRM's prediction is corrected and U1 heads on to $A_1$ 
following the adapted trajectory prediction. 
Please refer to Fig.\ref{fig:U1_hdg_to_A1_PR_offset_done_t8}.  \\
For the same reason, 
there is an PGRM process from $t_{140}$ to $t_{148}$, 
which corrects an offset track estimation for U1 leaving $A_1$ heading to $A_4$, 
see Fig.\ref{fig:U1_touchdown_t268}. \\
Exiting target area $A_4$, 
U1 follows the predicted trajectory as expected without 
violating the threshold $\kappa$. 
Subsequently, it visits $A_3$ and finally returns to $A_{base}$, 
which is reached at $t_{255}$. 
Successfully finishing mission execution, 
U1 adds a vertical landing phase and touches ground $@t_{268}$. 
Please note Fig.\ref{fig:U1_touchdown_t268}. \\
\subsection{Conclusions from the new scenario}
We proved that the Plan and Goal Recognition method presented in \cite{KreuzerWeberLeupolzKnoll23, KreuzerWeberKnoll24}, on which basis the PGRM is implemented, does not depend on an initial guess.\\
Additionally, we showed that even when trajectory candidates are near to identical, the impact of incorrect target identification on airspace management is negligible.
This is because short-term trajectory predictions remain accurate, naturally encompassing both potential tracks. \\
For medium-term intervals, any loss of validity ($d_{eukl.} > \kappa$) is immediately corrected through a triggered PGRM process.\\
As outcome, during the brief period in which the PGRM searches for a valid trajectory prediction, the system may experience a slight reduction in airspace safety.\\ 
%
%
%
%
\newpage
\bibliography{../../BibTex/mrabbrev,../../BibTex/articles,../../BibTex/books,../../BibTex/online}
\end{document}